\documentclass[draft]{publmathd}
\usepackage{amsmath,amsfonts,amssymb}

\pyear{2006} \pvol{1}

\refno{69/3}

\pagespan{291}{296}
\shortlefthead{A. Bovdi \es A. Szak\'acs}
\shortrighthead{Units of commutative group\dots}
\receivedinfo{January 4, 2006; revised September 1, 2006}

\newtheorem*{thm}{Theorem}
\newtheorem*{lem}{Lemma}

\newcommand\gp[1]{\langle#1\rangle}

\allowdisplaybreaks

\title[Units of commutative group algebra with involution]{Units of commutative group algebra with involution}

\author[A. Bovdi \es A. Szak\'acs]{A. \uppercase{Bovdi} (Debrecen) \es  A. \uppercase{Szak\'acs} (B\'ek\'escsaba)}

\address{Adalbert Bovdi\\
    Institute of Mathematics\\
    University of Debrecen\\
    H-4010 Debrecen, P.O. Box 12\\
    Hungary}

\email{bodibela@math.klte.hu}

\address{A. Szak\'acs\\
    Department of Business Mathematics\\
    Tessedik Samuel College\\
    H-5600 B\'ek\'escsaba\\
    Bajza u. 33.\\
    Hungary}

\email{szakacs@zeus.kf.hu}

\dedicatory {Dedicated to the memory of Dr.\ Edit Szab\'o}

\thanks{Supported by
OTKA No.\ T037202 and by FAPESP Brasil (proc.\ 06/56203-3)}
\subjclass{Primary: 16S34, 16U60; Secondary: 20C05} \keywords{group
algebra, group of units, unitary unit, symmetric unit}

\begin{document}

\begin{abstract}
Let  $p$ be an odd prime, $F$ the field of $p$ elements
and $G$ a finite abelian $p$-group with an arbitrary involutory
automorphism. Extend this automorphism to the group algebra $FG$
and consider the unitary and the symmetric normalized units of
$FG$. This paper provides bases and determines the invariants of
the two subgroups formed by these units.
\end{abstract}

\maketitle

\section{Introduction}

Let  $p$ be an odd prime, $F$ the field of $p$, $G$  a finite
abelian $p$-group with an arbitrary  automorphism $\eta$ of order
2 and $G_{\eta}=\{g\in G\mid \eta(g)=g\}$. Extending the
automorphism $\eta$ to the group algebra $FG$  we obtain the
involution
$$
x=\sum_{g\in G}\alpha_gg\mapsto x^{\circledast}=\sum_{g\in
G}\alpha_g\eta(g)
$$
of  $FG$  which we will call as $\eta$-cano\-nical involution. In
particular, if $\eta(g)=g^{-1}$ or $\eta(g)=g$  for all $g\in G$
then the involution is called involutory;  if $\eta(g)=g^{-1}$ for
all $g\in G$ then it is called cano\-nical and denote by $*$.

Let  $V(FG)=\{\sum_{g\in G}\alpha_gg\mid \sum_{g\in G}\alpha_g=1\}$ be the group of normalized units of the group
algebra $FG$ and  consider  the subgroups of symmetric and unitary units
\begin{align*}
S_{\circledast}(F G)&=\{ x \in V(F G) \mid x^{\circledast}=x\},\\[3pt]
 V_{\circledast}(F G)&=\{ x \in V(F G) \mid x^{\circledast}=x^{-1}\} \belowdisplayshortskip=-10pt \end{align*}
respectively.

Our goal is to study the unitary subgroup $V_{\circledast}(F G)$
and the group of symmetric units $S_{\circledast}(F G)$. The
problem of  determining the invariants and the basis of $V_*(F G)$
had been raised by S. P.~Novikov. Its solution for the canonical
involution was  given in  \cite{1}; here  this result extended to
arbitrary involutory involution.

\section{Invariants}

We start with some remarks about the invariants of unitary and symmetric subgroups to give a bases.

Since $V(FG)$ has an odd order and every $u\in V(F G)$  can be written as  $u=(v^{\circledast})^2$, so $
u=(v^{\circledast}v^{-1})(vv^{\circledast}), $ where $v^{\circledast}v^{-1}$ is unitary and $vv^{\circledast}$ is a
symmetric unit. But every  $x\in S_{\circledast}(F G) \cap V_{\circledast}(F G)$ is such that  $x=x^{\circledast}=x^{-1}$;
since $x$ is odd order it follows $x=1$ and we have
\begin{equation}
V(F G)=S_{\circledast}(F G)\times V_{\circledast}(F G). \end{equation}
Define   the mappings
$$
\psi_1: V(F G)\to V_{\circledast}(FG), \quad \psi_2: V(F G)\to S_{\circledast}(FG)
$$
given  by $\psi_1(x)=x^{\circledast} x^{-1}$ and
$\psi_2(x)=x^{\circledast} x$ respectively for $x\in V(F G)$. They
are epimorphisms and as corollary conclude  that
\begin{equation}
\begin{aligned}
V_{\circledast}(F G)&=\{x^{\circledast} x^{-1} \mid x \in V(FG)\},\\[3pt]
S_{\circledast}(F G)&=\{x^{\circledast} x \mid  x \in V(FG)\}\end{aligned} \end{equation}
and
\begin{equation}
S_{\circledast}(FG)^{p^i}=S_{\circledast}(F G^{p^i}),\quad  V_{\circledast}(FG)^{p^i}=V_{\circledast}(F G^{p^i}),\end{equation}
which use for the description of the invariants of these groups.
The subsets $\{g, \eta(g)\}$ with  $ g\in G\setminus G_\eta$ form
a partition of the set $G\setminus G_\eta$ and let $E$ be the
system of representatives of these subsets. Clearly,  $x\in S_{\circledast}(F G)$ can be uniquely  written as
$$
x=\sum_{g\in E}\alpha_g (g+\eta(g))+ \sum_{g\in G_\eta}\beta_gg $$
with  $\alpha_g, \beta_g\in F$ and $ \sum_{g\in E}2\alpha_g + \sum_{g\in G_\eta}\beta_g=1$, so  the order of the group of
symmetric units $S_{\circledast}(F G)$ equals  $ p^{\frac 1 2(|G|+|G_\eta|-2)}$. By (3)\quad  $S_{\circledast}(F G)^p=S_{\circledast}(F G^p)$,
so as before, the order of
$S_{\circledast}(F G)^p$ is $p^{ {1\over 2}(|G^p|+|G_\eta^p|-2)}$.
It follows that the $p$-rank of the group $S_{\circledast}(F G)$,
that is the number of components in the decomposition of
$S_{\circledast}(F G)$ into a direct product of cyclic groups,
equals $ \frac {1}{2}(|G|-|G^p|+|G_\eta|-|G_\eta^p|)$.  Similarly,
the $p$-rank of the group $S_{\circledast}(F G)^{p^{i-1}}$ equals
$ \frac {1}{2} (|G^{p^{i-1}}|- |G^{p^{i}}|+ |G_\eta^{p^{i-1}}|-|G_\eta^{p^{i}}|)$.

We conclude  that  the number of components of order $p^{i}$ in the decomposition of $S_{\circledast}(F G)$ into a
direct product of cyclic groups equals
\begin{equation}
f_i(S_{\circledast}(FG))= \frac{1}{2} (| G^{p^{i-1}}| -2 | G^{p^{i}} | + | G^{p^{i+1}}|+  | G_\eta^{p^{i-1}} |
-2 | G_\eta^{p^{i}} | + | G_\eta^{p^{i+1}}| ). \end{equation}
We known from \cite {3} that the equality $V(F G)^{p^{i}}= V(F G^{p^{i}})$ holds and the $p$-rank of $V(F G)^{p^{i-1}}$ equals
$|G^{p^{i-1}}|-|G^{p^{i}}|$. Now (1) it
yields that the $p$-rank of $V_{\circledast}(F G)^{p^{i-1}}$
equals $ \frac{1}{2} (|G^{p^{i-1}}|-|G^{p^{i}}|- |G_\eta^{p^{i-1}}|+|G_\eta^{p^{i}}| ). $ It  immediately follows
that  the number of components of order $p^{i}$ in the
decomposition of the group $V_{\circledast}(F G)$ into a direct product of cyclic groups is equal to
\begin{equation}
f_i(V_{\circledast}(FG))= \frac{1}{2}(| G^{p^{i-1}}| -2| G^{p^{i}} | + | G^{p^{i+1}} |- | G_\eta^{p^{i-1}}| +2 | G_\eta^{p^{i}}| - |
G_\eta^{p^{i+1}}|).\end{equation}

\section{The bases}

We will use the following well-known  generators of $V(FG)$ (see \cite{2}):

\begin{lem}
Let $G$ be a finite abelian $p$-group, $I=I(FG)$ the augmentation ideal of $FG$ and assume that $I^{s+1}=0$. If
$$
v_{d1}+I^{d+1},v_{d2}+I^{d+1},\dots,v_{d r_{d}}+I^{d+1} $$
is a basis for $I^{d}/I^{d+1}$, then the units $ 1+v_{dj}$ $(d=1,2,\dots, s$; $j=1,2,\dots,r_d) $ generate $V(F G)$.
\end{lem}

Let us return to the involutory automorphism $\eta$ of $G$.
Clearly,  $G$ has  the decomposition
$$
G=\langle a_1 \rangle\times \langle a_2 \rangle\times\dots \times\langle a_l\rangle\times\dots \times \langle a_t\rangle
$$
such that  the elements $a_1,a_2,\dots, a_l$ inverted by $\eta$
and if $t>l$ then $\eta$ leaves fixed $a_i$ for $i>l$.

Let $q_i$ be the order of $a_i$ and the set $L$ consisting of those
$t$-tuples $\alpha=(\alpha_1,\alpha_2,\dots ,\alpha_t)$ such that
$\alpha_i\in \{0,1,\dots,q_i-1 \}$ and at least one of $\alpha_j$ is
not divisible by $p$ ; the number of these elements $|G|-|G^p|$.
Write $L$ as the disjoint union $L=L_0\cup L_1\cup L_2$, where
$\alpha$ belongs to $L_0$, $L_1$ or $L_2$ according to whether
$\alpha_1+\alpha_2+\dots +\alpha_\ell$ is 0, odd or even and
positive. The cardinality of  $L_1$ is   $ {1\over
2}(|G|-|G^{p}|-|G_\eta|+|G_\eta^p |) $, and it is  the $p$-rank of
$V_{\circledast}(FG)$ and the cardinality of  $L_0\cup L_2$ is  $
 {1\over 2}(|G|-|G^{p}|+|G_\eta|-|G_\eta^p |)$, which is  the $p$-rank of $S_{\circledast}(FG)$.

\medskip
Put $ u_\alpha :=1+(a_{1}-1)^{\alpha_1}(a_{2}-1)^{\alpha_2} \dots (a_{t}-1)^{\alpha_t}$ for $\alpha\in L$.

\begin{thm}
Let $G$ be a finite abelian $p$-group of odd order with an involutory  automorphism  $\eta$ and $F$ is field of
$p$ elements. Then
\begin{itemize}
\item[{\rm 1.}]
The invariants of the unitary subgroup $V_{\circledast}(F G)$ are indicated in {\rm (5)} and the set
$\{{u_\alpha}^{\circledast} {u_\alpha}^{-1} \mid \alpha\in L_1 \}$ is basis for it, that is
$ V_{\circledast}(F G)=\prod_{\alpha \in L_1} \gp{{u_\alpha}^{\circledast} {u_\alpha}^{-1}}$.
\item[{\rm 2.}]
The invariants of the group of symmetric units $S_{\circledast}(FG )$ are indicated in {\rm (4)} and the set
$ \{u_\alpha^{\circledast} {u_\alpha} \mid \alpha\in L_2\}\cup \{u_\alpha \mid \alpha\in L_0\} $ is  basis for it, that is
$$
S_{\circledast}(F G)=\prod_{\alpha \in L_2} \gp{{u_\alpha}^{\circledast} {u_\alpha}}\times \prod_{\alpha \in L_0}  \gp{{u_\alpha}}.
$$
\end{itemize}
\end{thm}

\begin{proof}
Let $\alpha=(\alpha_1,  \alpha_2,  \dots, \alpha_t)$ with $\alpha_i\in \{0,1,\dots,q_i-1 \}$; define
$d(\alpha)={\alpha_1} + {\alpha_2} + \dots + {\alpha_t}$ and put
$ z_\alpha\,{:=}\,(a_{1}{-}\,1)^{\alpha_1}(a_{2}{-}\,1)^{\alpha_2} \dots (a_{t}{-}\,1)^{\alpha_t}$. It is known from \cite{2} that all
$z_\alpha$ form a vector space basis for the augmentation ideal
$I$ and the elements  $z_\alpha+I^{d+1}$ with $d(\alpha)=d$
constitute a basis for $I^d/I^{d+1}$ with the property required
for the application of Lemma. Therefore the elements
$u_\alpha=1+z_\alpha$ generate $V(FG)$ and from  \cite {3} follows
that  $\{u_\alpha | \alpha\in L\}$ is a basis of $V(FG)$.

\newpage
If  $i\leq l$ then  $(({a_i}-1) +1)^{-1}=(({a_i}-1)+1)^{\circledast}=(a_i-1)^{\circledast} +1$ and from
\begin{multline*}
(1+({a_i}-1))(1- (a_{i}-1)+({a_i}-1)^2 - \dots + (a_{i}-1)^{q_i-1})\\
=1+(a_{i}-1)^{q_i}=1 \end{multline*}
it follows that  $ (a_{i}-1)^{\circledast}=- (a_{i}-1)+({a_i}-1)^2 - \dots + (a_{i}-1)^{q_i-1}. $ Note that for  $i>l$  the equality
$(a_{i}-1)^{\circledast}=(a_{i}-1)$ holds.

Let $\alpha\in L_1\cup L_2$, $d={\alpha_1} + {\alpha_2} + \dots + {\alpha_t}$ and $k={\alpha_1} + {\alpha_2} + \dots + {\alpha_l}$.
The above argument ensures that $ {z_\alpha}^{\circledast} = (-1)^{k} z_{\alpha}+y$ for a  suitable  $y\in I^{d+1}$. It follows
that if $k$ is odd then
\begin{multline*}
(1+ z_{\alpha})^{-1}(1+ z_{\alpha})^{\circledast}=(1- z_{\alpha}+ z_{\alpha}^2-\dots)(1- z_{\alpha}+y)\\
\equiv 1-2 z_{\alpha}\pmod{I^{d+1}} \end{multline*} and for the
even $k$ we have
$$
(1+ z_{\alpha})(1+ z_{\alpha})^{\circledast}= (1+z_{\alpha}) (1+ z_{\alpha}+y)\equiv 1+2 z_{\alpha} \pmod{I^{d+1}}.
$$
Recall that $u_\alpha=1+z_\alpha$ and define
$$
z'_{\alpha}=\begin{cases} u_{\alpha}^{-1}u_{\alpha}^{\circledast}-1 , & \text{if} \ \alpha\in L_1;\\
u_{\alpha}u_{\alpha}^{\circledast}-1 , & \text{if} \ \alpha\in L_2;\\
z_\alpha,  &\text{if} \  \alpha\in L_0. \end{cases} $$
As a consequence of the foregoing  argument, we obtain,  modulo $I^{d+1}$,
$$
z'_{\alpha}\equiv \begin{cases} 2z_{\alpha} , &\text{if} \ \alpha\in L_1;\\
-2_{\alpha} , &\text{if} \ \alpha\in L_2;\\
z_\alpha,  &\text{if} \ \alpha\in L_0. \end{cases} $$
 Now  Lemma applies to the $z_\alpha'$, yielding  that the
 $u_\alpha'=1+z_\alpha'$ with $\alpha\in L$ also generate $V(FG)$.
 We claim that in fact they form a basis for $V(FG)$. To this end,
 it now suffices to show that the products of their orders
 is no larger than order $V(FG)$.
 Clearly that if $u_\alpha'\ne u_\alpha$ then such $u_\alpha'$ is the
 image of $u_\alpha$ either under the endomorphism
 $v\mapsto v^{\circledast}v^{-1}$ or under the endomorphism $v\mapsto  v^{\circledast} v$ of $V(FG)$ according to whether $\alpha$
 belongs to $L_1$ or $L_2$. It follows that  $|u_\alpha'|\leq |u_\alpha|$,
 where $|u_\alpha|$ is the order of $u_\alpha$.  But \cite{3}  asserts that $\prod_{\alpha\in L}|u_\alpha|$ is the order of $V(FG)$
and this gives that $\prod_{\alpha\in L}|u_\alpha'|\leq |V(FG)|$, so $\{u_\alpha' \mid \alpha\in L\}$ is a basis. It follows from
the definition of the $u_\alpha'$ with $\alpha\in L$ that each of
them is either fixed or inverted by the involution $\circledast$. Accordingly, this basis is the disjoint union of  bases for the
subgroups of symmetric and unitary normalized units, and these are
the bases in the theorem.
\end{proof}

\medskip
{\sc Acknowledgement.} The authors would like to thank the referee for their valuable comments and suggestions for clarifying the
exposition.

\end{document}